\documentclass[smallextended,numbook,runningheads]{svjour3}       % onecolumn (second format)
\smartqed  % flush right qed marks, e.g. at end of proof
\usepackage{graphicx}
\usepackage{amssymb,amsmath}
\usepackage{mathptmx}
\usepackage{graphics}
\newcommand{\dl}[1]{{\bf Theorem{#1}}}
\newcommand{\yl}[1]{{\bf Lemma{#1}}}

\newcommand{\de}[1]{{\bf Definition{#1}}}
\newcommand{\la}[1]{\label{#1}}
\newcommand{\rf}[1]{(\ref{#1})}
\newcommand{\ci}[1]{\cite{#1}}

\newcommand{\zm}{{\em Proof}}

\begin{document}

\title{The Gradient Superconvergence of Bilinear Finite Volume Element
for Elliptic Problems\thanks{This work was supported by the
National Natural Science Funds of China, No. 11371081; and the
State Key Laboratory of Synthetical Automation for Process
Industries Fundamental Research Funds, No. 2013ZCX02.}}

\titlerunning{Gradient superconvergence of bilinear finite
volume element}        % if too long for running head

\author{Tie Zhang \and Lixin Tang}
\authorrunning{T. Zhang, L.X. Tang} % if too long for running head
\institute{T. Zhang (Corresponding author) \at
              Department of Mathematics and the State Key Laboratory of Synthetical Automation for Process
Industries, Northeastern University, Shenyang 110004, China\\
              \email{ztmath@163.com}\qquad Tel.: +086-24-83680949\qquad
              Fax: +086-24-83680949
              \and
L.X. Tang \at  The State Key Laboratory of Synthetical Automation
for Process Industries,
              Northeastern University, Shenyang 110004,
China\\
\email{lixintang@mail.neu.edu.cn} }
\date{Received: date / Accepted: date}
\maketitle
\begin{abstract}
We study the gradient superconvergence of bilinear finite volume
element (FVE) solving the elliptic problems. First, a superclose
weak estimate is established for the bilinear form of the FVE
method. Then, we prove that the gradient approximation of the FVE
solution has the superconvergence property:
$\displaystyle{\max_{P\in S}}|(\nabla
u-\overline{\nabla}u_h)(P)|=O(h^2)|\ln h|$, where
$\overline{\nabla}u_h(P)$ denotes the average gradient on elements
containing point $P$ and $S$ is the set of optimal stress points
composed of the mesh points, the midpoints of edges and elements.

\end{abstract}

\section{Introduction}\label{intro}
The finite volume element (FVE) method has been widely used in
numerically solving partial deferential equations. The main
feature of FVE method is that it inherits some physical
conservation laws of original problems locally, which are very
desirable in practical applications. During the last decades, many
research works have been presented for FVE methods solving various
partial differential equations, see
\cite{Cai,Chen,Chen1,Chen2,Chen0,Chou,Chou1,Ewing,Lazarov,Li1,Lv,Suli,Wu,Xu,Zhang1,Zhang2,Zhang0}
and the references cited therein.

Superconvergence of numerical solutions has been an active
research area for finite element method (FEM) since its practical
importance in enhancing the accuracy of finite element
approximation. But the study of superconvergence properties of FVE
methods is far behind that of FEMs. For elliptic problems in
two-dimensional domain, the early superconvergence results of
linear and bilinear FVE solutions are of this form
\cite{Chen1,Li1}
\begin{equation}
\Big(\frac{1}{N}\sum_{z\in S}|(\nabla
u-\overline{\nabla}u_h)(z)|^2\Big)^{\frac{1}{2}}\leq
Ch^2\|u\|_{3,\infty},\label{1.1}
\end{equation}
where $S$ is the set of optimal stress points of interpolation
function on partition $T_h$, $N=O(h^{-2})$ is the total number of
points in $S$, and $\overline{\nabla}$ denotes the average
gradient on elements containing point $z$. Later, Lv and Li in
\cite{Lv} extended result \rf{1.1} to the isoparametric bilinear
FVE on quadrilateral meshes under the $h^2$-uniform mesh
condition. Recently, Zhang and Zou in \cite{Zhang0} also derived
some superconvergence results for the bi-complete $k$-order FVE on
rectangular meshes, and in the case of bilinear FVE ($k=1$), their
result is
\begin{equation}
|\nabla(u-u_h)(G)|\leq Ch^2|\ln
h|^{\frac{1}{2}}\|u\|_{4,\infty},\la{1.2}
\end{equation}
where $G$ is the Gauss point of element (the midpoint of element).
Moreover, by using the postprocessing technique, Chou and Ye
\cite{Chou1} obtain the superconvergence estimate:
\begin{equation}
\|\nabla u-\nabla Qu_h)\|\leq Ch^{\frac{3}{2}}\|u\|_3,\label{1.3}
\end{equation}
where $Qu_h$ is the postprocessed linear FVE solution obtained by
the $L_2$-projection method; Zhang and Sheng \cite{Zhang2} further
derive the superconvergence estimate:
\begin{equation}
\|\nabla u-R\nabla u_h)\|\leq Ch^2\|u\|_3,\label{1.4}
\end{equation}
where $R\nabla u_h$ is the reconstructed gradient of bilinear FVE
solution obtained by using the patch interpolation recovery
method.

In this paper, we consider the bilinear FVE method to solve the
following problem
\begin{eqnarray}
\left\{\begin{array}{ll}
-div(A\nabla u)+c\,u=f,\;in\;\;\Omega,\\
u=0,\;\;on\;\;\partial\Omega,\label{1.5}
\end{array}
 \right.
\end{eqnarray}
where $\Omega\subset R^2$ is a rectangular domain with boundary
$\partial\Omega$, coefficient matrix $A=(a_{ij})_{2\times 2}$. Our
main goal is to give some piecewise-point gradient
superconvergence for the bilinear FVE approximation to problem
\rf{1.5}. To the authors' best knowledge, in existing literatures,
only the midpoints of elements are proved to be the
superconvergence point of gradient approximation \cite{Zhang0}
(also see \rf{1.2}), we here will prove that except the midpoints
of elements, all interior mesh points and midpoints of interior
edges are also the superconvergence points. Generally speaking,
the analysis of bilinear FVE on rectangular meshes is more
difficult then that of linear FVE on triangle meshes, the reason
is that $\nabla u_h$ is not constant in the former case. By
calculating exactly some integrals on element and its boundary, we
first establish the superclose weak estimate for the bilinear form
of the FVE method,
\begin{equation}
|a_h(u-\Pi_hu,\Pi_h^*v)|\leq
Ch^2\|u\|_{3,p}\|v\|_{1,q},\;\forall\,v\in U_h,\;2\leq
p\leq\infty,\,\,1/p+1/q=1,\label{1.6}
\end{equation}
where $\Pi_h$ is the bilinear interpolation operator, $\Pi_h^*$ is
the projection operator from the trial function space $U_h$ to the
test function space. It is well known that such weak estimate
plays an important role in the superconvergence analysis of FEMs
\cite{Lin,Zhang,Zhu}. By means of this weak estimate and some
integral estimates on element, we further derive the following
superconvergence result for the gradient approximation
\begin{equation}
\max_{P\in S}|(\nabla u-\overline{\nabla}u_h)(P)|\leq Ch^2|\ln
h|\|u\|_{3,\infty},
\end{equation}
where $S$ is the optimal stress point set of bilinear
interpolation function which is composed of the interior mesh
points, the midpoints of interior edges and the midpoints of
elements.

This paper is organized as follows. In Section 2, we introduce the
FVE method and give some lemmas. In Section 3, we establish the
superconvergence weak estimate. Section 4 is devoted to the
piecewise-point gradient supconvergence estimate on the set $S$.
Finally, in Section 5, numerical experiments are provided to
illustrate our theoretical analysis.

Throughout this paper, we adopt the notations $W^{m,p}(D)$ to
indicate the usual Sobolev spaces on domain $D\subset \Omega$
equipped with the norm $\|\cdot\|_{m,p,D}$ and semi-norm
$|\cdot|_{m,p,D}$, and if $p=2$, we set $W^{m,p}(D)=H^m(D)$,
$\|\cdot\|_{m,p,D}=\|\cdot\|_{m,D}$. When $D=\Omega$, we omit the
index $D$. We will use letter $C$ to represent a generic positive
constant, independent of the mesh size $h$.
\section{Finite volume element method}
\setcounter{section}{2}\setcounter{equation}{0} Consider problem
\rf{1.5}. As usual, we assume that there exist positive constants
$C_1$ and $C_2$ such that
\begin{equation}
C_1\xi^T\xi\leq\xi^TA(x,y)\xi\leq C_2\xi^T\xi,\;\forall\,\xi\in
R^2,\;(x,y)\in\Omega,\label{2.1}
\end{equation}
We further assume that $A\in [W^{1,\infty}(\Omega)]^{2\times 2}$,
$c\in L_\infty(\Omega)$ and $c\geq 0$.

Let $T_h=\bigcup\lbrace K \rbrace$ be a rectangular partition of
domain $\Omega$ so that $\overline{\Omega}=\bigcup_{K\in
T_h}\lbrace\, \overline{K}\,\rbrace$, where $h=max \,h_K$, $h_K$
is the diameter of element $K$. We assume that partition $T_h$ is
regular, that is, there exists a positive constant $\gamma>0$ such
that
\begin{equation}
h_K/\rho_K\leq \gamma,\;\forall\, K\in T_h,\label{2.2}
\end{equation}
where $\rho_K$ denotes the diameter of the biggest ball included
in $K$.

Associated with partition $T_h$, we construct the central dual
partition $T^*_h$ by connecting the center of each element to the
midpoints of edges by straight lines. Thus, for each nodal point
$P$ in $T_h$, there exists a rectangle $K_P^*$ surrounding $P$,
$K_P^*\in T_h^*$ is called the dual element or the control volume
at point $P$, see Fig.1.\\
\parbox[c]{11cm}{
\begin{center}
\scalebox{0.5}{\includegraphics{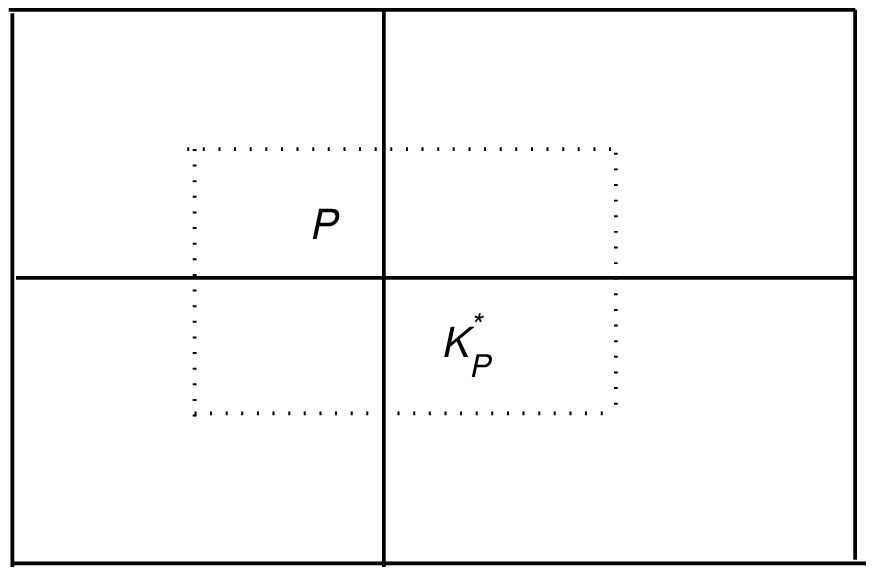}}
\end{center}

\centerline{\small  FIG.1.\quad The dual element $K_P^*$
surrounding point $P$}}

To partition $T_h$ and $T_h^*$, we associate the following trial
function space $U_h$ and test function space $V_h$, respectively,
\begin{eqnarray*}
&&U_h=\{\,u_h\in C^0(\overline{\Omega}):\,u_h|_K\in
Q_1(K),\;\forall\,K\in T_h,\;u_h|_{\partial\Omega}=0\,\},\\
&&V_h=\{\,v_h\in
L_2(\Omega):\,v_h|_{K_P^*}=constant,\,\forall\,P\in\,N_h,\;v_h|_{K_P^*}=0,\,\forall\,P\in\,\partial\Omega\},
\end{eqnarray*}
where $Q_1(K)$ is the set of all bilinear polynomials on $K$ and
$N_h$ is the set of all nodal points of $T_h$.

The conventional weak form for problem \rf{1.5} is that find $u\in
H^1_0(\Omega)$ such that
\begin{equation}
a(u,v)=(f,v),\;\forall\,v\in H^1_0(\Omega),\label{2.3}
\end{equation}
where
\begin{equation}
a(u,v)=\int_{\Omega}A\nabla u\cdot\nabla
v+cu\,v,\;\;(f,v)=\int_{\Omega}f\,v\,\label{2.4}.
\end{equation}
Weak form \rf{2.3} usually is adopted for the finite element
method. But for the FVE method, we need a new weak form. Let $u$
be the solution of problem \rf{1.5}, then by using Green's
formula, we have
\begin{equation}
-\int_{\partial K_P^*}n\cdot(A\nabla
u)vds+\int_{K_P^*}cuv=\int_{K_P^*}fv,\;\;\,K_P^*\in T^*_h,\;v\in
V_h,\label{2.5}
\end{equation}
where $n$ is the outward unit normal vector on the boundary
concerned. Motivated by this weak formula, we introduce the
following bilinear form
\begin{equation}
a_h(u,v)=\sum_{K_P^*\in T^*_h}\Big(-\int_{\partial K_P^*}n\cdot
(A\nabla u)vds+\int_{K_P^*}cuv\Big),\;u\in\,H^1(\Omega),\,v\in
V_h,\label{2.6}
\end{equation}
and define the FVE approximation of problem \rf{1.5} by finding
$u_h\in U_h$ such that
\begin{equation}
a_h(u_h,v_h)=(f,v_h),\;\forall\,v_h\in V_h.\label{2.7}
\end{equation}

Let $\Pi_h^*:\,U_h\rightarrow V_h$ be the interpolation operator
defined by
$$
\Pi_h^*v_h=\sum_{P\in N_h}v_h(P)\chi_{_P},\;\forall\,v_h\in U_h,
$$
where $\chi_{_P}$ is the characteristic function of the dual
element $K_P^*$. Since $\Pi_h^*$ is a one to one mapping from
$U_h$ onto $V_h$, then the equivalent form of problem \rf{2.7} is
that find $u_h\in U_h$ such that
\begin{equation}
a_h(u_h,\Pi_h^*v_h)=(f,\Pi_h^*v_h),\;\forall\,v_h\in
U_h.\label{2.8}
\end{equation}
This is the FVE scheme to be used in our analysis. From \rf{2.5}
we know that scheme \rf{2.8} is consistent and the following error
equation holds.
\begin{equation}
a_h(u-u_h,\Pi_h^*v_h)=0,\;\forall\, v_h\in U_h.\label{2.9}
\end{equation}

Let $\Pi_hu$ be the usual bilinear interpolation function of
continuous function $u$. In our analysis, the following
approximation property, trace inequality and inverse inequality
will be used frequently. For $1\leq p,q\leq\infty$, we have
\begin{eqnarray}
&&\|u-\Pi_hu\|_{m,p,K}\leq
Ch_K^{2-m}\|u\|_{2,p,K},\;0\leq m\leq 2,\;1<p\leq\infty,\;\label{2.10}\\
&&\|u\|_{0,p,\partial K}\leq
Ch_K^{-\frac{1}{p}}\big(\,\|u\|_{0,p,K}+h_K\|\nabla
u\|_{0,p,K}\big),\;u\in W^{1,p}(K),\label{2.11}\\
&&\|u_h\|_{m,p,K}\leq
Ch_K^{\frac{2}{p}-\frac{2}{q}}h_K^{l-m}\|u_h\|_{l,q,K},\;u_h\in
Q_1(K),\;0\leq l\leq m\leq 1 .\label{2.12}
\end{eqnarray}
Furthermore, the following
two lemmas hold \cite{Li1,Zhang2}.\\
\yl{ 2.1}\quad{\em Let $K\in T_h$, $\tau\subset\partial K$ be an
edge of $K$. Then, for $v_h\in U_h,\,1\leq q\leq \infty$, we have}
\begin{eqnarray}
&&\int_K(v_h-\Pi_h^*v_h)=0,\;\;\int_\tau(v_h-\Pi_h^*v_h)ds=0,\label{2.13}\\
&&\|v_h-\Pi_h^*v_h\|_{0,q,K}\leq Ch_K\|v_h\|_{1,q,K},\label{2.14}\\
&&\|v_h-\Pi_h^*v_h\|_{0,q,\partial K}\leq
Ch_K^{1-\frac{1}{q}}\|v_h\|_{1,q,K}.\label{2.15}
\end{eqnarray}
\yl{ 2.2}{\em \quad For $h$ small, we have}
$$
a_h(v_h,\Pi_h^*v_h)\geq C\|v_h\|^2_1,\;\forall\,v_h\in U_h.
$$

Lemma 2.2 implies that the solution $u_h$ of problem \rf{2.8}
uniquely exists.

\section{Interpolation weak estimate}
\setcounter{section}{3}\setcounter{equation}{0} It is well known
that the interpolation weak estimate plays an important role in
the superconvergence analysis \cite{Lin,Zhang,Zhu}. To FEM defined
on regular rectangular meshes, the following interpolation weak
estimate has been established \cite{Zhang,Zhu}
\begin{equation}
|a(u-\Pi_hu,v)|\leq Ch^2\|u\|_{3,p}\|v\|_{1,q},\;\forall\,v\in
U_h,\;2\leq p\leq\infty,\,1/p+1/q=1\,.\label{3.1}
\end{equation}
In this section, we will give a similar weak estimate for the FVE
method, which is very useful in our superconvergence analysis.

In order to utilize the known result \rf{3.1}, we need
to give the difference between the bilinear form of FVE and that of FEM.\\
\yl{ 3.1}\quad{\em For any $w\in U_h+H^2(\Omega),\,v_h\in U_h$, we
have
\begin{eqnarray}
a_h(w,\Pi_h^*v_h)-a(w,v_h)&=&\sum_{K\in T_h}\int_{\partial
K}n\cdot(A\nabla
w)(\Pi_h^*v_h-v_h)ds\nonumber\\
&+&\sum_{K\in T_h}(-div(A\nabla w)+cw
,\Pi_h^*v_h-v_h)_K,\label{3.2}
\end{eqnarray}
where $U_h+H^2(\Omega)=\{ w=u_h+v: u_h\in U_h, v\in H^2(\Omega)
\}$ is the the algebraic sum space.}\\
\zm\quad By Green's formula, we have
\begin{eqnarray*}
\int_KA\nabla w\cdot\nabla v_h=-\int_Kdiv(A\nabla
w)v_h+\int_{\partial K}n\cdot(A\nabla w)v_hds,
\end{eqnarray*}
and (see Fig.1)
\begin{eqnarray*}
&&\sum_{K\in T_h}\int_Kdiv(A\nabla w)\Pi_h^*v_h=\sum_{K\in
T_h}\sum_{K_P^*\in T_h^*}\int_{K_P^*\cap K}div(A\nabla w)\Pi_h^*v_h\\
&=&\sum_{K\in T_h}\int_{\partial K}n\cdot(A\nabla w)\Pi_h^*v_hds
+\sum_{K_P^*\in T^*_h}\int_{\partial K_P^*}n\cdot(A\nabla
w)\Pi_h^*v_hds.
\end{eqnarray*}
Substituting this two identities into the definitions of
$a(w,v_h)$ and $a_h(w,\Pi_h^*v_h)$ (see \rf{2.4} and \rf{2.6}),
the proof is completed.$\hfill\square$

Let $K=\Box P_1P_2P_3P_4$ be a rectangular element,
$h_x=x_2-x_1,\,h_y=y_2-y_1$, see Fig.2. The following lemma gives
some exact calculations for integrals on $\partial K$.

\parbox[c]{12cm}{
\begin{center}
\scalebox{0.5}{\includegraphics{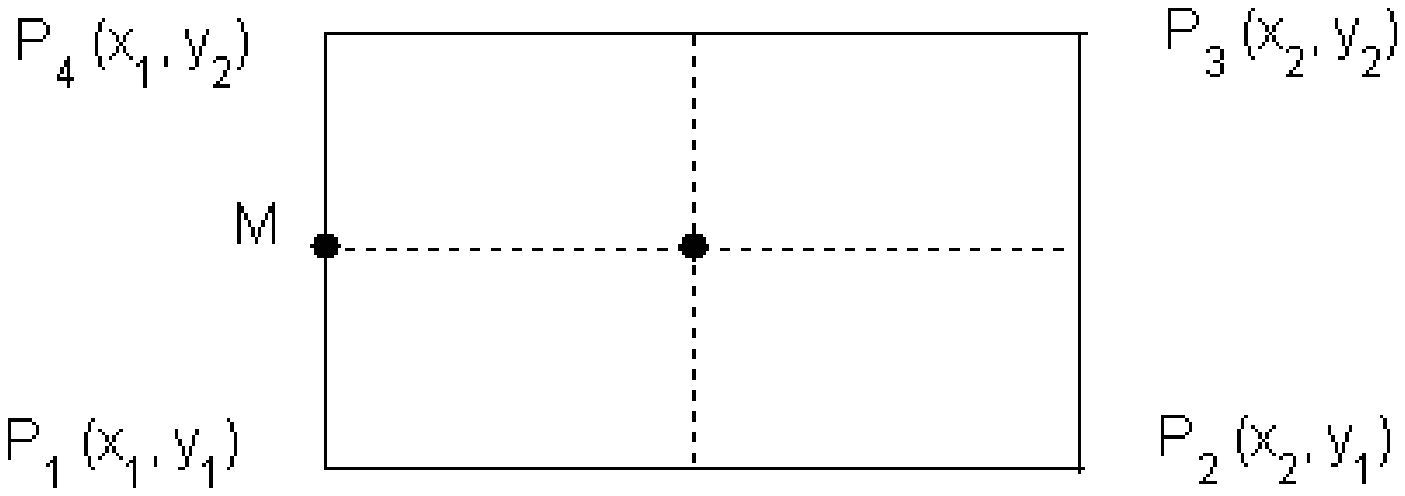}}
\end{center}
} \centerline{\small  FIG.2.\quad Rectangular element $K=\Box
P_1P_2P_3P_4$}\\
\yl{ 3.2}\quad{\em  For $w,v\in U_h$, we have}
\begin{eqnarray}
\int_{P_1P_4}(\Pi_h^*v-v)w_xdy=\frac{h_y^3}{24}v_y(x_1,y)w_{xy},\;\;\int_{P_1P_4}(\Pi_h^*v-v)w_ydy=0,\label{3.3}\\
\int_{P_2P_3}(\Pi_h^*v-v)w_xdy=\frac{h_y^3}{24}v_y(x_2,y)w_{xy},\;\;\int_{P_2P_3}(\Pi_h^*v-v)w_ydy=0,\label{3.4}\\
\int_{P_1P_2}(\Pi_h^*v-v)w_ydx=\frac{h_x^3}{24}v_x(x,y_1)w_{xy},\;\;\int_{P_1P_2}(\Pi_h^*v-v)w_xdx=0,\label{3.5}\\
\int_{P_4P_3}(\Pi_h^*v-v)w_ydx=\frac{h_x^3}{24}v_x(x,y_2)w_{xy},\;\;\int_{P_4P_3}(\Pi_h^*v-v)w_xdx=0.\label{3.6}
\end{eqnarray}
\zm\quad We only prove \rf{3.3}, the other are similar. Noting
that $w_x$ and $v$ are linear on line segment $P_1P_4$ ($M$ is its
middle point), we have from the definition of $\Pi_h^*v$ that
\begin{eqnarray*}
\int_{P_1P_4}\Pi_h^*vw_xdy&=&\int_{P_1M}v(P_1)w_xdy+\int_{MP_4}v(P_4)w_xdy\\
&=&\frac{h_y}{4}v(P_1)\big(w_x(P_1)+w_x(M)\big)+\frac{h_y}{4}v(P_4)\big(w_x(M)+w_x(P_4)\big)\\
&=&\frac{h_y}{4}\big(v(P_1)w_x(P_1)+2v(M)w_x(M)+v(P_4)w_x(P_4)\big).
\end{eqnarray*}
On the other hand, using the composite trapezoidal rule:
$$
\int_a^bf(y)dy=\frac{b-a}{4}\big(f(a)+2f(\frac{a+b}{2})+f(b)\big)-\frac{(b-a)^3}{48}f''_y
$$
we obtain by taking $f=vw_x$ that
\begin{eqnarray*}
\int_{P_1P_4}vw_xdy=\frac{h_y}{4}\big(v(P_1)w_x(P_1)+2v(M)w_x(M)+v(P_4)w_x(P_4)\big)-\frac{h_y^3}{24}v_y(x_1,y)w_{xy}.
\end{eqnarray*}
The difference of above two equalities implies the first equality
in \rf{3.3}. The second equality comes from \rf{2.13}, noting that
$w_y$ is constant on $P_1P_4$.$\hfill\square$

Let $w^c$ denote the piecewise constant approximation of function
$w$ on $T_h$,
\begin{equation}
w^c|_K=\frac{1}{K}\int_Kw,\;K\in T_h;\;\;\|w-w^c\|_{0,p,K}\leq
Ch_K\|w\|_{1,p,K},\;1\leq p\leq\infty.\label{3.7}
\end{equation}
\dl{ 3.1}{\em \quad Let $T_h$ be a regular rectangular mesh and
$u\in W^{3,p}(\Omega)$. Then we have}
\begin{equation}
|a_h(u-\Pi_hu,\Pi_h^*v)|\leq
Ch^2\|u\|_{3,p}\|v\|_{1,q},\;\forall\,v\in U_h,\;2\leq
p\leq\infty,\,1/p+1/q=1\,.\label{3.8}
\end{equation}
\zm\quad Let $A_M=(\overline{a}_{ij})_{2\times 2}$ be the value of
$A$ at the midpoint $M$ of edge $\tau\subset\partial K$, see
Fig.2. From Lemma 3.1, we have
\begin{eqnarray}
&&a_h(u-\Pi_hu,\Pi_h^*v)-a(u-\Pi_hu,v)\nonumber\\
&=&\sum_{K\in T_h}\int_{\partial K}n\cdot(A-A_M)\nabla
(u-\Pi_hu)(\Pi_h^*v-v)ds\nonumber\\
&&+\sum_{K\in T_h}\int_{\partial
K}n\cdot(A_M\nabla (u-\Pi_hu))(\Pi_h^*v-v)ds\nonumber\\
&&+\sum_{K\in T_h}(-div(A\nabla
(u-\Pi_hu)),\Pi_h^*v-v)_K\nonumber\\
&&+\sum_{K\in T_h}(c(u-\Pi_hu),\Pi_h^*v-v)_K\nonumber\\
&=&E_1+E_2+E_3+E_4\label{3.9}.
\end{eqnarray}
Using \rf{2.10}-\rf{2.11} and \rf{2.14}-\rf{2.15}, we obtain
\begin{eqnarray*}
E_1&\leq& C\sum_{K\in
T_h}h_K|A|_{1,\infty}\|\nabla(u-\Pi_hu)\|_{0,p,\partial
K}\|v-\Pi_h^*v\|_{0,q,\partial K}\leq
Ch^2\|u\|_{2,p}\|v\|_{1,q}.\\
E_4&\leq& C\sum_{K\in
T_h}\|u-\Pi_hu\|_{1,p,K}\|v-\Pi_h^*v\|_{0,q,K}\leq
Ch^2\|u\|_{2,p}\|v\|_{1,q}.
\end{eqnarray*}
For $E_3$, set $A=({\bf a_1,a_2})$. Since
\begin{equation}
div(A\nabla w)=(div\,{\bf a_1}, div\,{\bf a_2})\cdot\nabla w+{\bf
a_1}\cdot\nabla w_x+{\bf a_2}\cdot\nabla w_y,\label{3.10}
\end{equation}
then we have
\begin{eqnarray*}
E_3&=&\sum_{K\in T_h}(-(div\,{\bf a_1}, div\,{\bf a_2})\cdot\nabla
(u-\Pi_hu),\Pi_h^*v-v)_K\\
&&-\sum_{K\in T_h}({\bf a_1}\cdot\nabla(u-\Pi_hu)_x+{\bf
a_2}\cdot\nabla (u-\Pi_hu)_y,\Pi_h^*v-v)_K=E_{31}+E_{32}.
\end{eqnarray*}
Obviously, $E_{31}\leq Ch^2\|u\|_{2,p}\|v\|_{1,q}$. Now, let
$u_i=u(P_i),\,\xi=(x-x_1)/h_x,\;\eta=(y-y_1)/h_y$, see Fig.2.
Since
$$
\Pi_hu=u_1+(u_2-u_1)\xi+(u_4-u_1)\eta+(u_3+u_1-u_2-u_4)\xi\eta,\;(x,y)\in
K,
$$
hence
\begin{eqnarray}
(\Pi_hu)_{xy}&=&(u_3+u_1-u_2-u_4)h_x^{-1}h_y^{-1}\nonumber\\
&=&\frac{1}{K}
\big(\int_{x_1}^{x_2}u_x(x,y_2)dx-\int_{x_1}^{x_2}u_x(x,y_1)dx\big)\nonumber\\
&=&\frac{1}{K}
\int_{x_1}^{x_2}\int_{y_1}^{y_2}u_{xy}dxdy=\frac{1}{K}\int_Ku_{xy}dxdy=(u_{xy})^c.\label{3.11}
\end{eqnarray}
Therefore, we have from \rf{2.13} that
\begin{eqnarray*}
E_{32}&=&-\sum_{K\in
T_h}(a_{11}u_{xx}+a_{21}(u_{xy}-(\Pi_hu)_{xy})+
a_{12}(u_{xy}-(\Pi_hu)_{xy})+a_{22}u_{yy},\Pi_h^*v-v)_K\\
&=&-\sum_{K\in T_h}(a_{11}u_{xx}+a_{21}(u_{xy}-u_{xy}^c)+ a_{12}
(u_{xy}-u_{xy}^c)+a_{22}u_{yy},\Pi_h^*v-v)_K\\
&=& -\sum_{K\in
T_h}(a_{11}u_{xx}+a_{22}u_{yy}-(a_{11}u_{xx}+a_{22}u_{yy})^c,\Pi_h^*v-v)_K\\
&&-\sum_{K\in T_h}(a_{21}(u_{xy}-u_{xy}^c)+ a_{12}
(u_{xy}-u_{xy}^c),\Pi_h^*v-v)_K.
\end{eqnarray*}
Using the approximation property, we obtain $E_{32}\leq
Ch^2\|u\|_{3,p}\|v\|_{1,q}$ so that
$$
E_{3}=E_{31}+E_{32}\leq Ch^2\|u\|_{3,p}\|v\|_{1,q}.
$$
Now we need to estimate
\begin{eqnarray}
E_2=\sum_{K\in T_h}\sum_{\tau\subset\partial
K\setminus\partial\Omega}\int_\tau n\cdot(A_M\nabla
(u-\Pi_hu))(\Pi_h^*v-v)ds,\label{3.12}
\end{eqnarray}
noting that $(\Pi_h^*v-v)|_{\partial\Omega}=0$. Let $\tau$ be an
interior edge, that is, a common edge of two adjacent element $K$
and $K'$. Since $n|_{\tau\cap\partial K}=-n|_{\tau\cap\partial
K'}$ and $A_M\nabla u(\Pi_h^*v-v)$ is continuous across edge
$\tau$, then we have
\begin{equation}
E_2=\sum_{K\in T_h}\sum_{\tau\subset\partial
K\setminus\partial\Omega}\int_\tau -n\cdot(A_M\nabla
\Pi_hu)(\Pi_h^*v-v)ds\equiv\sum_{K\in
T_h}\sum_{\tau\subset\partial
K\setminus\partial\Omega}F(\tau).\label{3.13}
\end{equation}
Without loss of generality, let $K=\Box
P_1P_2P_3P_4=(x_1,x_2)\times(y_1,y_2)$ be a rectangular element
and the common edge $\tau=\partial K\bigcap\partial K'=P_1P_4$,
see Fig.2. From Lemma 3.2 we have
\begin{equation}
F(\tau)=\int_{P_1P_4}(\overline{a}_{11}(\Pi_hu)_x+\overline{a}_{12}(\Pi_hu)_y)
(\Pi_h^*v-v)dy=\overline{a}_{11}\frac{h_y^3}{24}v_y(x_1,y)(\Pi_hu)_{xy}.\label{3.14}
\end{equation}
Then, from \rf{3.14}, \rf{3.4} and \rf{3.11}, and noting that
$n|_{\tau\cap\partial K}=-n|_{\tau\cap\partial K'}$,
$\overline{a}_{11}v_y(x_1,y)$ is continuous acrose $\tau=\partial
K\cap\partial K'=P_1P_4$, we obtain
\begin{eqnarray}
F(\tau\cap\partial K)+F(\tau\cap\partial
K')=\overline{a}_{11}\frac{h_y^3}{24}v_y(x_1,y)\big((u_{xy})^c(K)-(u_{xy})^c(K')\big).\label{3.15}
\end{eqnarray}
Let $i_hw\in C^0(\Omega),\,i_hw|_K\in Q_1(K)$, be an approximation
of function $w$, which satisfies
$$
\|w-i_hw\|_{0,p}\leq Ch\|w\|_{1,p},\;2\leq p\leq\infty.
$$
For example, for $p=2$ we may take $i_hw=P_hw$, $P_h$ is the
$L_2$-projection operator; For $p>2$, we may take $i_hw=\Pi_hw$.
Then, from \rf{3.15} and the inverse inequality, we have
\begin{eqnarray}
&&F(\tau\cap\partial K)+F(\tau\cap\partial
K')\nonumber\\
&&=\overline{a}_{11}\frac{h_y^3}{24}v_y(x_1,y)\big((u_{xy})^c(K)-i_hu_{xy}(M)+i_hu_{xy}(M)-(u_{xy})^c(K')\big)
\nonumber\\
&&\leq
Ch_y^3h_K^{-\frac{2}{q}}\|v\|_{1,q,K}\big(h_K^{-\frac{2}{p}}\|(u_{xy})^c-i_hu_{xy}\|_{0,p,K}
+h_{K'}^{-\frac{2}{p}}
\|(u_{xy})^c-i_hu_{xy}\|_{0,p,K'}\big)\nonumber\\
&&\leq Ch_y\|v\|_{1,q,K}\big(\|(u_{xy})^c-u_{xy}\|_{0,p,K\cup K'
}+\|u_{xy}-i_hu_{xy}\|_{0,p,K\cup K'}\big),\label{3.16}
\end{eqnarray}
where we have used the fact that $h_y=|\partial K\bigcap\partial
K'|\leq \min\{h_K,h_{K'}\}$. Combining \rf{3.13} with \rf{3.16},
we can conclude that
\begin{equation}
E_2\leq
Ch\|v\|_{1,q}\big(\|u_{xy}-(u_{xy})^c\|_{0,p}+\|u_{xy}-i_hu_{xy}\|_{0,p}\big)\leq
Ch^2\|u\|_{3,p}\|v\|_{1,q}.\label{3.17}
\end{equation}
Substituting estimates $E_1\sim E_4$ into \rf{3.9}, we complete
the proof by using \rf{3.1}.$\hfill\square$

From Theorem 3.1, we immediately obtain the following superclose
result.\\
\dl{ 3.2}{\em \quad Let $u$ and $u_h$ be the solutions of problem
\rf{1.5} and \rf{2.8}, respectively, $u\in H^3(\Omega)$. Then we
have }
$$
\|\Pi_hu-u_h\|_1\leq Ch^2\|u\|_3.
$$
\zm\quad From Lemma 2.2, error equation \rf{2.9} and weak estimate
\rf{3.8}, we have
\begin{eqnarray*}
&&C\|u_h-\Pi_hu\|_1^2\leq a_h(u_h-\Pi_hu,\Pi_h^*(u_h-\Pi_hu))\\
&&=a_h(u-\Pi_hu,\Pi_h^*(u_h-\Pi_hu))\leq
Ch^2\|u\|_3\|u_h-\Pi_hu\|_1.
\end{eqnarray*}
This gives the conclusion.\hfill$\square$

By using the embedding theory and weak embedding inequality of
finite element space \ci{Zhu}, a direct result from Theorem 3.2 is
the following optimal $L_p$-error estimates,
\begin{eqnarray*}
&&\|u-u_h\|_{0,p}\leq \|u-\Pi_hu\|_{0,p}+\|\Pi_hu-u_h\|_{0,p}\leq
C_ph^2(\|u\|_{2,p}+\|u\|_3),\,1<p<\infty,\\
&&\|u-u_h\|_{0,\infty}\leq
C\big(h^2\|u\|_{2,\infty}+|ln\,h|^{\frac{1}{2}}\|\Pi_hu-u_h\|_1\big)\leq
Ch^2(\|u\|_{2,\infty}+|ln\,h|^{\frac{1}{2}}\|u\|_3).
\end{eqnarray*}

\section{Piecewise-point superconvergence for gradient approximation}
\setcounter{section}{4}\setcounter{equation}{0} \de{ 4.1}\quad{\em
Let point $P\in \Omega$, we call $P$ an optimal stress point for
interpolation operator $\Pi_h$, if
\begin{equation}
|(\nabla u-\overline{\nabla}\Pi_hu)(P)|\leq
Ch^2\|u\|_{3,\infty,E},\label{4.1}
\end{equation}
where $E$ is the union of elements containing point $P$ and
$\overline{\nabla}w(P)$ denotes the average value of gradient
$\nabla w(P)$ in elements containing point $P$.}

It is well known that for the bilinear interpolation $\Pi_hu$ on
rectangular meshes, the optimal stress point set $S$ is composed
of all interior mesh points, midpoints of interior edges and the
midpoints of elements \cite{Zhang3,Zhu}. In this section, we will
derive the $W^{1,\infty}$-superconvergence result for the bilinear
finite volume solution $u_h$ on set $S$.

We first give some lemmas. Let rectangular element $K=\square
P_1P_2P_3P_4$. For $w_h\in U_h$, set $w_i=w(P_i), i=1,\cdots,4$.
Introduce the discrete $H^1$-norm by
\begin{equation}
|w_h|_{1,h}=\big(\sum_{K\in
T_h}|w_h|_{1,h,K}^2\big)^{\frac{1}{2}},\;w_h\in U_h,\label{4.2}
\end{equation}
where
\begin{equation}
|w_h|_{1,h,K}^2=(w_2-w_1)^2+(w_3-w_2)^2+(w_3-w_4)^2+(w_4-w_1)^2.\label{4.3}
\end{equation}
\yl{ 4.1}\quad{\em The discrete norm $|w_h|_{1,h}$ is equivalent
to the norm $|w_h|_1$ on space $U_h$, and}
\begin{equation}
\frac{1}{6\gamma}|w_h|_{1,h,K}^2\leq |w_h|_{1,K}^2\leq
\frac{\gamma}{2}|w_h|^2_{1,h,K},\;\forall\, K\in T_h.\label{4.4}
\end{equation}
\zm\quad Let $K=\square P_1P_2P_3P_4=[x_1,x_1+h_x]\times
[y_1,y_1+h_y], \xi=(x-x_1)/h_x,\;\eta=(y-y_1)/h_y$. For $w_h\in
Q_1(K)$, we may write
\begin{equation}
w_h(x,y)=w_1+(w_2-w_1)\xi+(w_4-w_1)\eta+(w_3+w_1-w_2-w_4)\xi\eta.\label{4.5}
\end{equation}
Set
\begin{eqnarray}
&&w_{21}=w_2-w_1,\,w_{41}=w_4-w_1,\,w_{32}=w_3-w_2,\,w_{34}=w_3-w_4,\label{4.6}\\
&&w_{1234}=w_3+w_1-w_2-w_4=w_{34}-w_{21}=w_{32}-w_{41}.\label{4.7}
\end{eqnarray}
Then, we have from \rf{4.5}-\rf{4.7} that
\begin{eqnarray*}
|w_h|_{1,K}^2&=&\int_K|\nabla
w_h|^2dxdy=\int_0^1\int_0^1\big((\partial_\xi
w_h/h_x)^2+(\partial_\eta w_h/h_y)^2\big)h_xh_yd\xi d\eta\\
&=&\frac{h_y}{h_x}\int_0^1\int_0^1(w_{21}+w_{1234}\eta)^2d\xi
d\eta+
\frac{h_x}{h_y}\int_0^1\int_0^1(w_{41}+w_{1234}\xi)^2d\xi d\eta\\
&=&\frac{h_y}{h_x}(w_{21}^2+w_{21}w_{1234}+\frac{1}{3}w_{1234}^2)
+\frac{h_x}{h_y}(w_{41}^2+w_{41}w_{1234}+\frac{1}{3}w_{1234}^2)\\
&=&\frac{1}{3}\frac{h_y}{h_x}(w_{21}^2+w_{21}w_{34}+w_{34}^2)
+\frac{1}{3}\frac{h_x}{h_y}(w_{32}^2+w_{32}w_{41}+w_{41}^2).
\end{eqnarray*}
Hence, it follows from the Cauchy inequality and the regularity
condition \rf{2.2} that
$$
\frac{1}{6\gamma}(w_{21}^2+w_{34}^2+w_{32}^2+w_{41}^2)\leq
|w_h|_{1,K}^2\leq
\frac{\gamma}{2}(w_{21}^2+w_{34}^2+w_{32}^2+w_{41}^2).
$$
Combining this with \rf{4.3}, the proof is
completed.\hfill$\square$

Let $A^c=(a^c_{ij})_{2\times 2}$ be the piecewise constant
approximation of matrix function $A$.\\
\yl{ 4.2}\quad{\em The following identity holds for $w,v\in U_h$},
\begin{eqnarray}
\int_{\partial K}n\cdot (A^c\nabla
w)(\Pi_h^*v-v)ds=\frac{h_y^3h_x}{24}a^c_{11}v_{xy}w_{xy}+\frac{h_yh_x^3}{24}a^c_{22}v_{xy}w_{xy}.\label{4.8}
\end{eqnarray}
\zm\quad From Lemma 3.2 we obtain
\begin{eqnarray*}
&&\int_{\partial K}n\cdot (A^c\nabla w)(\Pi_h^*v-v)ds\\
&=&\int_{P_2P_3}(\Pi_h^*v-v)(a^c_{11}w_x+a^c_{12}w_y)dy-\int_{P_1P_4}(\Pi_h^*v-v)(a^c_{11}w_x+a^c_{12}w_y)dy\\
&&+\int_{P_4P_3}(\Pi_h^*v-v)(a^c_{21}w_x+a^c_{22}w_y)dx-\int_{P_1P_2}(\Pi_h^*v-v)(a^c_{21}w_x+a^c_{22}w_y)dx\\
&=&\frac{h_y^3}{24}a^c_{11}v_y(x_2,y)w_{xy}-\frac{h_y^3}{24}a^c_{11}v_y(x_1,y)w_{xy}\\
&&+\frac{h_x^3}{24}a^c_{22}v_x(x,y_2)w_{xy}-\frac{h_x^3}{24}a^c_{22}v_x(x,y_1)w_{xy}\\
&=&\frac{h_y^3h_x}{24}a^c_{11}v_{xy}w_{xy}+\frac{h_yh_x^3}{24}a^c_{22}v_{xy}w_{xy}.\
\hbox{\hspace{6.2cm}$\square$}
\end{eqnarray*}

Now, we need to introduce the regularized Green function
\cite{Zhang,Zhu}. For any given $z\in \Omega$, let $\delta_h^z\in
U_h$ be the smooth $\delta$-function which satisfies,
$$
(\delta_h^z,v_h)=v_h(z),\;z\in\Omega,\;\forall\, v_h\in U_h.
$$
For any appointed direction $L$, define the direction derivative
$$
\partial_zv(z)=\lim_{\vartriangle z\rightarrow 0,\vartriangle
z//L}\big(v(z+\vartriangle z)-v(z)\big)/|\vartriangle z|.
$$
Then, there exists a regularized Green function of derivative type
$\partial_z G^z(x)\in H^1_0(\Omega)\bigcap H^2(\Omega)$ such that
$$
a(v,\partial_zG^z)=(\partial_z\delta^z_h,v),\;\forall\,v\in
H^1_0(\Omega).
$$
Let $\partial_zG_h^z\in U_h$ be the finite element approximation
of $\partial_zG^z$ such that
$$
a(v_h,\partial_zG^z-\partial_zG_h^z)=0,\;\forall\,v_h\in U_h.
$$
Clearly, we have
\begin{equation}
a(v_h,\partial_zG^z_h)=a(v_h,\partial_zG^z)=(\partial_z\delta^z_h,v_h)=\partial_zv_h(z),
\;\forall\,v_h\in U_h.\label{4.9}
\end{equation}
Assume that partition $T_h$ is quasi-uniform so that the inverse
inequality holds on finite element space $U_h$. Under this
condition, the following boundness estimates were given in
\cite{Zhang,Zhu}
\begin{equation}
\|\partial_z G^z_h\|_1\leq Ch^{-1}|\ln
h|^{\frac{1}{2}},\;\;\;\|\partial_z G^z_h\|_{1,1}\leq C|\ln
h|,\label{4.10}
\end{equation}
where $C$ is a positive constant independent of $z\in\Omega$.\\
\dl{ 4.1}{\em \quad Let partition $T_h$ be quasi-uniform, and $u$
and $u_h$ be the solutions of problem \rf{1.5} and \rf{2.8},
respectively, $u\in W^{3,\infty}(\Omega)$. Then we have}
\begin{equation}
\|\Pi_hu-u_h\|_{1,\infty}\leq Ch^2|\ln
h|\|u\|_{3,\infty}.\label{4.11}
\end{equation}
\zm\quad From \rf{4.9}, error equation \rf{2.9}, Theorem 3.1 and
\rf{4.10}, we obtain
\begin{eqnarray}
&&\partial_z(u_h-\Pi_hu)(z)=a(u_h-\Pi_hu,\partial_z G^z_h)\nonumber\\
&&=a(u_h-\Pi_hu,\partial_z G^z_h)-a_h(u_h-\Pi_hu,\Pi_h^*\partial_z
G^z_h)+a_h(u_h-\Pi_hu,\Pi_h^*\partial_z G^z_h)\nonumber\\
&&=a(u_h-\Pi_hu,\partial_z G^z_h)-a_h(u_h-\Pi_hu,\Pi_h^*\partial_z
G^z_h)+a_h(u-\Pi_hu,\Pi_h^*\partial_z G^z_h)\nonumber\\
&&\leq
a(u_h-\Pi_hu,\partial_z G^z_h)-a_h(u_h-\Pi_hu,\Pi_h^*\partial_z
G^z_h)+Ch^2\|u\|_{3,\infty}\|\partial_z G^z_h\|_{1,1}\nonumber\\
&&\leq a(u_h-\Pi_hu,\partial_z
G^z_h)-a_h(u_h-\Pi_hu,\Pi_h^*\partial_z G^z_h)+Ch^2|\ln
h|\|u\|_{3,\infty}.\label{4.12}
\end{eqnarray}
Below we need to estimate (see Lemma 3.1)
\begin{eqnarray}
&&a(u_h-\Pi_hu,\partial_z G^z_h)-a_h(u_h-\Pi_hu,\Pi_h^*\partial_z
G^z_h)\nonumber\\
&=&\sum_{K\in T_h}\int_{\partial K}n\cdot(A-A^c)\nabla
(u_h-\Pi_hu)(\partial_z G^z_h-\Pi_h^*\partial_z
G^z_h)ds\nonumber\\
&&+\sum_{K\in T_h}\int_{\partial K}n\cdot A^c\nabla
(u_h-\Pi_hu)(\partial_z G^z_h-\Pi_h^*\partial_z
G^z_h)ds\nonumber\\
&&+\sum_{K\in T_h}(-div(A\nabla (u_h-\Pi_hu)),\partial_z
G^z_h-\Pi_h^*\partial_z
G^z_h)_K\nonumber\\
&&+(c(u_h-\Pi_hu),\partial_z G^z_h-\Pi_h^*\partial_z
G^z_h)_K\nonumber\\
&=&S_1+S_2+S_3+S_4.\label{4.13}
\end{eqnarray}
First, using \rf{2.11}, \rf{2.12} and \rf{2.15}, Theorem 3.2 and
\rf{4.10}, we have
\begin{eqnarray*}
S_1&\leq& C\sum_{K\in
T_h}h_K|A|_{1,\infty}\|\nabla(u_h-\Pi_hu)\|_{0,\partial
K}\|\partial_z G^z_h-\Pi_h^*\partial_z G^z_h\|_{0,\partial K}\\
&\leq&
Ch^{\frac{1}{2}}\|\nabla(u_h-\Pi_hu)\|h^{\frac{1}{2}}\|\partial_z
G^z_h\|_1\\
&\leq&Ch^3\|u\|_3\|\partial_z G^z_h\|_1\leq Ch^2|\ln
h|^{\frac{1}{2}}\|u\|_3.
\end{eqnarray*}
Next, it follows from Lemma 4.2 and inverse inequality \rf{2.12}
that
\begin{eqnarray*}
S_2&=&\leq \sum_{K\in
T_h}|A|_{\infty}\frac{h_K^4}{12}|(u_h-\Pi_hu)_{xy}|_{0,\infty,K}|(\partial_z
G^z_h)_{xy}|_{0,\infty,K}\\
&\leq& C\sum_{K\in
T_h}h_K^2|(u_h-\Pi_hu)_{xy}|_{0,\infty,K}\|\partial_z
G^z_h\|_{1,K}.
\end{eqnarray*}
From \rf{4.5} and Lemma 4.1, we have
$$
|(u_h-\Pi_hu)_{xy}|=|(u_h-\Pi_hu)_{1234}|\leq
\sqrt{2}\,|u_h-\Pi_hu|_{1,K,h}\leq
\sqrt{12\gamma}\;|u_h-\Pi_hu|_{1,K}.
$$
Substituting this estimate into $S_2$ and using Theorem 3.2, we
obtain
$$
S_2\leq Ch^2\|u_h-\Pi_hu\|_1\|\partial_z G^z_h\|_1\leq Ch^3|\ln
h|^{\frac{1}{2}}\|u\|_3.
$$
Similarly, by Lemma 4.1 and Theorem 3.2, we have (noting that
$(u_h-\Pi_hu)_{xx}=(u_h-\Pi_hu)_{yy}=0$\,)
\begin{eqnarray*}
S_3&\leq& C\sum_{K\in
T_h}\big(\|\nabla(u_h-\Pi_hu)\|_{0,K}+\|(u_h-\Pi_hu)_{xy}\|_{0,K}\big)\|\partial_z
G^z_h-\Pi_h^*\partial_z G^z_h)\|_{0,K}\\
&\leq& C\sum_{K\in T_h}\|u_h-\Pi_hu\|_{1,K}\|\partial_z
G^z_h-\Pi_h^*\partial_z G^z_h)\|_{0,K}\\
&\leq &Ch^3\|u\|_3\|\partial_z G^z_h\|_1\leq Ch^2|\ln
h|^{\frac{1}{2}}\|u\|_3.
\end{eqnarray*}
For $S_4$, using Theorem 3.2, \rf{2.14} and \rf{4.10} to obtain
\begin{eqnarray*}
S_4&\leq& C\sum_{K\in
T_h}|A|_{1,\infty}\|u_h-\Pi_hu\|_{1,K}\|\partial_z
G^z_h-\Pi_h^*\partial_z G^z_h\|_{0,K}\\
&\leq& Ch^2\|u\|_3\,h\|\partial_z G^z_h\|_1\leq Ch^2|\ln
h|^{\frac{1}{2}}\|u\|_3.
\end{eqnarray*}
Substituting estimates $S_1\sim S_4$ into \rf{4.13}, we complete
the proof by \rf{4.12}.\hfill$\square$

Now we can give the main superconvergence result.\\
\dl{ 4.2}{\em \quad Let partition $T_h$ be quasi-uniform, and $u$
and $u_h$ be the solutions of problem \rf{1.5} and \rf{2.8},
respectively, $u\in W^{3,\infty}(\Omega)$. Then we have
\begin{equation}
\max_{P\in S}|(\nabla u-\overline{\nabla}u_h)(P)|\leq Ch^2|\ln
h|\|u\|_{3,\infty},\label{4.14}
\end{equation}
where $\overline{\nabla}w(P)$ denotes the average values of
gradient $\nabla w(P)$ in elements containing point $P$ and $S$ is
the optima stress point set composed of all interior mesh points,
the midpoints of interior edges and the midpoints of elements.}\\
\zm\quad From \rf{4.1} and Theorem 4.1, we immediately obtain the
conclusion by using the triangle inequality.\hfill$\square$

\section{Numerical example}
\setcounter{section}{5}\setcounter{equation}{0}
In this section, we
will present some numerical results to illustrate the theoretical
analysis.

Let us consider problem \rf{1.5} with the data:
$$
A(x,y)=\left(
\begin{array}{cc}
  e^{2x}+y^3+1 & e^{x+y} \\
  e^{x+y} & e^{2y}+x^3+1
\end{array}
\right),\;\;\;\; c(x,y)=2+x+y.
$$
We take $\Omega=[0,1]^2$ and the exact solution $u(x,y)=2\sin(2\pi
x)\sin(3\pi y)$.

In the numerical experiments, we first partition domain $\Omega$
into a rectangular mesh $T_h$, then the refined meshes are
obtained by connecting the midpoints of each edge of elements in
$T_h$. Thus, a rectangular mesh sequence is generated with
successively halving mesh size $h/2^i,\,i=1,2,\cdots$.

Denote by $e_h=\displaystyle{\max_{P\in S}}|\nabla
u(P)-\overline{\nabla }u_h(P)|$ the computation error and the
numerical convergence rate is computed by using the formula
$r=\ln(e_h/e_{h/2})/\ln 2$. Table 1 gives the error and
convergence rate with successively halved mesh sizes. We see that
the convergence rate is just about $O(h^2)$, as
the theoretical prediction.\\
%\vspace{0.2cm}
\begin{center}
{\small {\bf Table 1}\quad Convergence rate of gradient approximation on set $S$}\\[0.5\baselineskip]
\renewcommand\arraystretch{1.2}
\arrayrulewidth 0.5pt
\begin{tabular}{ccccccc}
\hline
$h$& 1/4& 1/8 & 1/16 & 1/32 & 1/64 &1/128\\
\hline
$e_h$&1.212&3.099e-1&7.856e-2&1.969e-2&4.949e-3&1.243e-3\\

rate &--&1.9671&1.9802&1.9961&1.9926&1.9932\\
\hline
\end{tabular}
\end{center}
%\vskip 0.2cm

\end{document}